\documentclass[10pt]{article}

\usepackage{a4wide}
\usepackage{amssymb}
\usepackage{amsfonts}
\usepackage{amsmath}
\input xy
\xyoption{arrow} \xyoption{matrix}
\usepackage[all,cmtip]{xy}

\date{}

\newtheorem{proposition}{Proposition}[section]
\newtheorem{theorem}[proposition]{Theorem}
\newtheorem{lemma}[proposition]{Lemma}

\newtheorem{corollary}[proposition]{Corollary}

\def\der{\partial }

\def\nFM0{{\nu }_{F,M_0}}
\def\nFN0{{\nu }_{F,N_0}}
\def\nGN0{{\nu }_{G,N_0}}

\def\N0{ {\bf N}_0 }

\def\ra{\rightarrow}

\def\Xpm{X^{\pm }}

\def\s{\sigma}

\def\l1{{\lambda}_1}

\def\a{\alpha}
\def\a0{ {\alpha }_0}
\def\a1{ {\alpha }_1}

\def\l{\lambda}

\def\nFGM0{{\nu }_{F,G,M_0}}

\def\nFN0{{\nu}_{F,N_0}}

\def\sm{{\sigma}^m}

\def\sm1{{\sigma}^{-1}}

\def\smtp1{{\sigma}^{-t+1}}

\def\S1{S^{-1}}

\def\Xpm1{X^{\pm 1}_1}

\def\sPM1{{\sigma }^{\pm 1}}
\def\sMP1{{\sigma }^{\mp 1 }}

\def\d{\delta}

\def\di{{\rm d.ind}}

\def\L{\Lambda}

\def\Ytm1{Y^{t-1}}
\def\Yim1{Y^{i-1}}

\def\CL{{\cal L}}

\def\CN{{\cal N}}

\def\ass{{\rm ass}}

\def\ker{ {\rm ker } }

\def\gr{ {\rm gr} }

\def\SL2Z{ {\rm SL}_2({\bf Z}) }

\def\th{ \theta }

\def\CL{{\cal L}}

\def\Gp1{ G^{1 , 1 } }
\def\P11{ P^{-1 , 1 } }
\def\Pp1{ P^{1 , 1 } }

\def\th{\theta}

\def\nCLsr{{}^\nu\kern-2pt {\cal L}^{\sigma , \rho  }}
\def\nP{{}^\nu \kern-2pt P}
\def\nL{{}^\nu\kern-2pt L}
\def\nLL{{}^\nu\kern-2pt \Lambda}
\def\nPsr{{}^\nu\kern-2pt P^{\sigma , \rho  }}
\def\nLsr{{}^\nu\kern-2pt L^{\sigma , \rho  }}
\def\nuCL{{}^\nu\kern-2pt  {\cal L}}
\def\nCLsr{{}^\nu\kern-2pt {\cal L}^{\sigma , \rho  }}
\def\nCL1m{{}^\nu\kern-2pt {\cal L}^{-1 , 1  }}
\def\x1nu{x^\frac{1}{\nu}}
\def\xm1nu{x^{-\frac{1}{\nu}}}

\def\rad{{\rm rad}}

\def\CN{{\cal N}}
\def\ra{\rightarrow }

\def\CB{{\cal B}}

\def\CT{{\cal T}}

\def\CC{ {\cal C}}

\def\nAM0{{\nu }_{{\cal A},M_0}}
\def\nAN0{{\nu }_{{\cal A},N_0}}

\def\bR{\overline{R}}

\def\ga{\mathfrak{a}}
\def\gb{\mathfrak{b}}
\def\gc{\mathfrak{c}}

\def\gn{\mathfrak{n}}

\def\gr{\mathfrak{r}}

\def\SL{{\rm SL}}

\def\di!{\frac{\der^i}{i!}}
\def\dik!{\frac{\der^k_i}{k!}}

\def\N{\mathbb{N}}

\def\0{\overline{0}}
\def\1{\overline{1}}

\def\Ln1{\L_{n,\overline{1}}}

\def\a1{a_{\overline{1}}}

\def\S{\Sigma}

\def\vn1{\overrightarrow{n-1}}

\def\im{{\rm im}}

\def\Min{{\rm Min}}

\def\mJ{\mathbb{J}}
\def\mI{\mathbb{I}}

\def\K1{{\rm K}_1}

\def\hmI1{\widehat{\mI_1}}
\def\tmI1{\widetilde{\mI_1}}
\def\tmJ1{\widetilde{\mJ_1}}
\def\hB1{\widehat{B_1}}
\def\hCB1{\widehat{\CB_1}}

\def\Den{{\rm Den}}
\def\Denl{{\rm Den}_l}
\def\Ore{{\rm Ore}}

\def\Den{{\rm Den}}

\def\Loc{{\rm Loc}}

\def\Ass{{\rm Ass}}

\def\maxDen{{\rm max.Den}}
\def\maxAss{{\rm max.Ass}}
\def\maxLoc{{\rm max.Loc}}
\def\llrad{{\rm l.lrad}}

\def\assmaxDen{{\rm ass.max.Den}}

\def\br{\overline{r}}

\def\ga{\mathfrak{a}}

\def\udim{{\rm udim}}

\def\bE{\overline{E}}

\def\Nil{{\rm Nil}}
\def\gll{\mathfrak{l}}

\begin{document}

\author{V. V. \  Bavula 
}

\title{Weakly left localizable rings}

\maketitle

\begin{abstract}

A new class of rings, {\em the class of weakly  left localizable rings},  is
introduced. A ring $R$ is called  {\em  weakly left localizable} if each non-nilpotent  element
of $R$ is invertible in some left localization $S^{-1}R$ of the
ring $R$. Explicit criteria are given for  a ring to be a weakly left
localizable ring provided the ring has only finitely many maximal
left denominator sets (eg, this is the case if a ring has a left
Artinian left quotient ring). It is proved that a ring with finitely
many maximal left denominator sets that satisfies some natural conditions  is a weakly left localizable ring iff
 its left quotient
ring is a direct product of finitely many  local  rings such that their radicals are nil ideals.




$${\bf Contents}$$
\begin{enumerate}
\item Introduction.
\item Preliminaries.
\item Weakly left localizable rings and their characterizations.
 \item The cores of maximal denominator sets of a weakly left localizable ring.
 \item Criterion for a semilocal ring to be a  weakly left localizable ring.
\end{enumerate}
\end{abstract}


\section{Introduction}


Throughout, module means a left module. In this paper the following notation will remained fixed.

$\noindent $

{\bf Notation}:

\begin{itemize}

\item $R$ is a ring with 1, $R^*$ is its  group of units  and
$\CC_R$ is the set of (left and right) regular elements of the
ring $R$  (i.e. $\CC_R$ is the set of non-zero-divisors of $R$);
\item $\rad (R)$ is the radical of $R$,   $\CN_R$ is the nil radical of $R$ and $\Nil (R)$ is the set of nilpotent elements of $R$;
\item  $Q_{l, cl}(R):=\CC_R^{-1}R$ is the {\em classical left quotient ring} of $R$ (if it exists);
\item $\Ore_l(R):=\{ S\, | \, S$ is a left Ore set in $R\}$; \item
$\Den_l(R):=\{ S\, | \, S$ is a left denominator set in $R\}$;
\item $\Loc_l(R):= \{ [S^{-1}R]\, | \, S\in \Den_l(R)\}$ where $[S^{-1}R]$ is an $R$-{\em isomorphism} class of the ring $S^{-1}R$ (a ring isomorphism $\s : S^{-1}R\ra S'^{-1}R$ is called an $R$-{\em isomorphism} if $\s (\frac{r}{1})= \frac{r}{1}$ for all elements $r\in R$);
    \item
$\Ass_l(R):= \{ \ass (S)\, | \, S\in \Den_l(R)\}$ where $\ass
(S):= \{ r\in R \, | \, sr=0$ for some $s=s(r)\in S\}$; \item
$\Den_l(R, \ga ) := \{ S\in \Den_l(R)\, | \, \ass (S)=\ga \}$
where $\ga \in \Ass_l(R)$;
\item $\maxDen_l(R)$ is the set of maximal left denominator sets of $R$;
\item $\gll_R:= \llrad (R):= \bigcap_{S\in \maxDen_l(R)}\ass (S)$ is the {\em left localization radical} of the ring $R$;
\item      $\Denl
(R,0)$ is the set of regular  left denominator sets $S$ in $R$
($S\subseteq \CC_R$);
\item $S_\ga=S_\ga (R)=S_{l,\ga }(R)$
is the {\em largest element} of the poset $(\Den_l(R, \ga ),
\subseteq )$ and $Q_\ga (R):=Q_{l,\ga }(R):=S_\ga^{-1} R$ is  the
{\em largest left quotient ring associated to} $\ga$, $S_\ga $
exists (Theorem 
 2.1, \cite{larglquot});
\item In particular, $S_0=S_0(R)=S_{l,0}(R)$ is the largest
element of the poset $(\Den_l(R, 0), \subseteq )$ and
$Q_l(R):=S_0^{-1}R$ is the largest left quotient ring of $R$;
\end{itemize}


{\bf Weakly left localization maximal rings and their characterization}.  The aim of the paper is to introduce a new class of rings, the class of weakly left localizable rings, and to give several  characterizations of them (Theorem \ref{26Mar14} and Theorem \ref{28Mar14}) in the case when they admit only finitely many maximal left denominator sets and satisfy some natural conditions. Notice that each ring with left Artinian, left quotient ring has only finitely many maximal left denominator sets, \cite{Bav-LocArtRing}. Let $R$ be a ring, an element $r\in R$ is called a {\em left localizable element},   \cite{Crit-S-Simp-lQuot}, if  there exists a left denominator set $S= S(r)$ such that $r\in S$ (equivalently, there exists a left denominator set $S'= S'(r)$ such that $\frac{r}{1}$ is a unit in $S'^{-1}R$). Let $\CL_l(R)$ be the {\em set of left localizable elements} of the ring $R$ and $\CN\CL_l(R):= R\backslash \CL_l(R)$ be the {\em set of left non-localizable elements} of $R$. A ring $R$ is called a {\em left localizable ring} if $\CL_l(R)=R\backslash \{ 0\}$, \cite{left-localizable-rings}. In \cite{left-localizable-rings}, several characterizations of the class of left
 localizable rings are given.

 \begin{itemize}
\item ({\bf Theorem 3.9},  \cite{left-localizable-rings}) {\em Let $R$ be a ring. The following statements are equivalent.}
\begin{enumerate}
\item {\em The ring $R$ is a left localizable ring with
$n:=|\maxDen_l(R)|<\infty$.\item $Q_{l,cl}(R)=R_1\times\cdots
\times R_n$  where $R_i$ are division rings. }
 \item {\em The ring $R$ is a semiprime left Goldie ring with $\udim
 (R)=|\Min (R)|=n$ where $\Min (R)$ is the set of minimal prime
 ideals of the ring $R$.}
\item $Q_l(R)= R_1\times\cdots \times R_n$ {\em  where $R_i$  are division rings.}
\end{enumerate}
\end{itemize}

A ring $R$ is called  a {\em local ring} if the set $R\backslash R^*$ is an ideal of the ring $R$ (equivalently, $R/ \rad (R)$ is a division ring).
\begin{itemize}
\item ({\bf Theorem \ref{26Mar14}})
{\em Let $R$ be a ring. The following statements are equivalent.}
\begin{enumerate}
\item {\em The ring $R$ is a weakly left localizable ring such that}
\begin{enumerate}
\item $\gll_R=0$,
\item $|\maxDen_l(R)|<\infty$,
\item {\em for every $S\in \maxDen_l(R)$, $S^{-1}R$ is a weakly left localizable ring, and
    \item for all $S, T\in \maxDen_l(R)$ such that $S\neq T$, $\ass (S)$ is not a nil ideal  modulo $\ass (T)$}.
\end{enumerate}
\item $Q_{l, cl}(R)\simeq \prod_{i=1}^n R_i$ {\em  where $R_i$ are local rings with $\rad (R_i) = \CN_{R_i}$.
\item $Q_l(R)\simeq \prod_{i=1}^n R_i$ where $R_i$ are local rings with } $\rad (R_i) = \CN_{R_i}$.
\end{enumerate}

\end{itemize}

\begin{itemize}
\item ({\bf Theorem \ref{28Mar14}})  {\em Let $R$ be a ring, $\gll = \gll_R$, $\pi': R\ra R':= R/\gll$, $r\mapsto \br := r+\gll$.
The following statements are equivalent.}
\begin{enumerate}
\item {\em The ring $R$ is a weakly left localizable ring such that}
\begin{enumerate}
\item {\em the map $\phi : \maxDen_l(R)\ra  \maxDen_l(R')$, $S\mapsto \pi' (S)$, is a surjection.
\item $|\maxDen_l(R)|<\infty$,
\item for every $S\in \maxDen_l(R)$, $S^{-1}R$ is a weakly left localizable ring, and
    \item  for all $S, T\in \maxDen_l(R)$ such that $S\neq T$, $\ass (S)$ is not a nil ideal  modulo $\ass (T)$}.
\end{enumerate}
\item $Q_{l, cl}(R')\simeq \prod_{i=1}^n R_i$ {\em where $R_i$ are local rings with $\rad (R_i) = \CN_{R_i}$, $\gll$ is a nil ideal and $\pi'(\CL_l(R))= \CL_l(R')$.
\item $Q_l(R')\simeq \prod_{i=1}^n R_i$ where $R_i$ are local rings with $\rad (R_i) = \CN_{R_i}$,  $\gll$ is a nil ideal and } $\pi'(\CL_l(R))= \CL_l(R')$.
\end{enumerate}
\end{itemize}

Rings that satisfy the assumptions of Theorem \ref{26Mar14} or Theorem \ref{28Mar14} have interesting properties (see Corollary \ref{b26Mar14}, Corollary \ref{c28Mar14} and Corollary \ref{d28Mar14}).

\begin{itemize}
\item ({\bf Corollary  \ref{b26Mar14}}) {\em  We keep the notation of Theorem \ref{26Mar14}. Suppose that a ring $R$ satisfies one of the equivalent conditions 1--3 of Theorem \ref{26Mar14}. Then}
\begin{enumerate}
\item $\maxDen_l(R)=\{ S_1, \ldots , S_n\}$ {\em where $S_i=\{ r\in R\, | \, \frac{r}{1}\in R_i^*\}$.
\item $\CC_R=\bigcap_{S\in \maxDen_l(R)}S$.
\item $\Nil (R) = \CN_R$.
\item $Q:= Q_{l,cl} (R) = Q_l(R)$ is a weakly left localizable ring with} $\Nil (Q) = \CN_Q= \rad (Q)$.
    \item $\CC_R^{-1}\CN_R =\CN_Q
     = \rad (Q)$.
    \item $\CC_R^{-1} \CL_l(R) = \CL_l(Q)$.
\end{enumerate}
\end{itemize}

{\bf The cores of maximal left denominator sets of weakly left localizable rings}.  Let $R$ be a ring and $S$ be its left Ore set. The subset of $S$,
$$ S_c:=\{ s\in S\, | \, \ker (s\cdot ) = \ass (S)\}$$ is called the {\em core} of the left Ore set $S$ where $s\cdot : R\ra R$, $r\mapsto sr$, \cite{Crit-S-Simp-lQuot}. In \cite{left-localizable-rings}, properties of the core of a left Ore set are studied in detail.


\begin{itemize}
\item ({\bf Theorem 4.2},  \cite{left-localizable-rings})
{\em Suppose that $S\in \Den_l(R, \ga )$ and $S_c\neq \emptyset$. Then}
\begin{enumerate}
\item $S_c\in \Den_l(R, \ga )$.
\item {\em The map $\th : S_c^{-1}R\ra S^{-1}R$, $s^{-1}r\mapsto s^{-1}r$, is a ring isomorphism. So,} $S_c^{-1}R\simeq S^{-1}R$.
\end{enumerate}
\end{itemize}

The next theorem gives an explicit description of the cores of maximal left denominator sets of a weakly left localizable ring $R$ that satisfies the conditions of Theorem \ref{26Mar14}.
\begin{itemize}
\item ({\bf Theorem \ref{C2Dec12}})
{\em Let $R$ be a ring  that satisfies one of the equivalent  conditions 1--3 of Theorem \ref{26Mar14}. Let } $\maxDen_l(R)=\{ S_1, \ldots , S_n\}$.
\begin{enumerate}
\item {\em If $n=1$ then $S_{1,c} = S_1= R\backslash \CN_R$.
\item If $n\geq 2$ then $S_{i, c} = S_i\cap \bigcap_{j\neq i} \ga_i\neq \emptyset $  where $\ga_j=\ass (S_j)$.}
\end{enumerate}
\end{itemize}
{\bf Criterion for a semilocal ring to be a weakly left localizable ring}.
A ring $R$ is called a {\em semilocal ring} if $R/\rad (R)$ is a semisimple (Artinian) ring. The next theorem is a criterion for a semilocal ring $R$ to be a weakly left localizable ring with $\rad (R) = \CN_R$.

\begin{itemize}
\item ({\bf Theorem \ref{24Dec12}})
{\em Let $R$ be a semilocal ring. Then the ring $R$ is a weakly left localizable ring with $\rad (R)=\CN_R$  iff $R\simeq \prod_{i=1}^s R_i$ where $R_i$ are local rings with} $\rad (R_i) =\CN_{R_i}$.
\end{itemize}


\section{Preliminaries}\label{PRREML}

 In this section,  necessary results are collected  that are used in the proofs of the paper.


{\bf The largest regular left Ore set and the largest left
quotient ring of a ring}. Let $R$ be a ring. A {\em
multiplicatively closed subset} $S$ of $R$ or a {\em
 multiplicative subset} of $R$ (i.e. a multiplicative sub-semigroup of $(R,
\cdot )$ such that $1\in S$ and $0\not\in S$) is said to be a {\em
left Ore set} if it satisfies the {\em left Ore condition}: for
each $r\in R$ and
 $s\in S$,
$$ Sr\bigcap Rs\neq \emptyset .$$
Let $\Ore_l(R)$ be the set of all left Ore sets of $R$.
  For  $S\in \Ore_l(R)$, $\ass (S) :=\{ r\in
R\, | \, sr=0 \;\; {\rm for\;  some}\;\; s\in S\}$  is an ideal of
the ring $R$.


A left Ore set $S$ is called a {\em left denominator set} of the
ring $R$ if $rs=0$ for some elements $ r\in R$ and $s\in S$ implies
$tr=0$ for some element $t\in S$, i.e. $r\in \ass (S)$. Let
$\Den_l(R)$ be the set of all left denominator sets of $R$. For
$S\in \Den_l(R)$, let $S^{-1}R=\{ s^{-1}r\, | \, s\in S, r\in R\}$
be the {\em left localization} of the ring $R$ at $S$ (the {\em
left quotient ring} of $R$ at $S$).

In general, the set $\CC$ of regular elements of a ring $R$ is
neither left nor right Ore set of the ring $R$ and as a
 result neither left nor right classical  quotient ring ($Q_{l,cl}(R):=\CC^{-1}R$ and
 $Q_{r,cl}(R):=R\CC^{-1}$) exists.
 Remarkably, there  exists the largest
 regular left Ore set $S_0= S_{l,0} = S_{l,0}(R)$. This means that the set $S_{l,0}(R)$ is an Ore set of
 the ring $R$ that consists
 of regular elements (i.e., $S_{l,0}(R)\subseteq \CC$) and contains all the left Ore sets in $R$ that consist of
 regular elements. Also, there exists the largest regular (left and right) Ore set $S_{l,r,0}(R)$ of the ring $R$.
 In general, all the sets $\CC$, $S_{l,0}(R)$, $S_{r,0}(R)$ and $S_{l,r,0}(R)$ are distinct, for example,
 when $R= \mI_1:=K\langle x, \frac{d}{dx}, \int\rangle$ is the ring of polynomial integro-differential operators over a field $K$ of characteristic zero,
  \cite{intdifline}.

$\noindent $

{\it Definition}, \cite{intdifline}, \cite{larglquot}.    The ring
$$Q_l(R):= S_{l,0}(R)^{-1}R$$ (respectively, $Q_r(R):=RS_{r,0}(R)^{-1}$ and
$Q(R):= S_{l,r,0}(R)^{-1}R\simeq RS_{l,r,0}(R)^{-1}$) is  called
the {\em largest left} (respectively, {\em right and two-sided})
{\em quotient ring} of the ring $R$.

$\noindent $

 In general, the rings $Q_l(R)$, $Q_r(R)$ and $Q(R)$
are not isomorphic, for example, when $R= \mI_1$, \cite{intdifline}.  The next
theorem gives various properties of the ring $Q_l(R)$. In
particular, it describes its group of units.


\begin{theorem}\label{4Jul10}
\cite{larglquot}
\begin{enumerate}
\item $ S_0 (Q_l(R))= Q_l(R)^*$ {\em and} $S_0(Q_l(R))\cap R=
S_0(R)$.
 \item $Q_l(R)^*= \langle S_0(R), S_0(R)^{-1}\rangle$, {\em i.e. the
 group of units of the ring $Q_l(R)$ is generated by the sets
 $S_0(R)$ and} $S_0(R)^{-1}:= \{ s^{-1} \, | \, s\in S_0(R)\}$.
 \item $Q_l(R)^* = \{ s^{-1}t\, | \, s,t\in S_0(R)\}$.
 \item $Q_l(Q_l(R))=Q_l(R)$.
\end{enumerate}
\end{theorem}

The set $(\Den_l(R), \subseteq )$ is a poset (partially ordered
set). In \cite{larglquot}, it is proved  that the set
$\maxDen_l(R)$ of its maximal elements is a {\em non-empty} set.

{\bf The maximal denominator sets and the maximal left localizations  of a ring}.


{\it Definition}, \cite{larglquot}. An element $S$ of the set
$\maxDen_l(R)$ is called a {\em maximal left denominator set} of
the ring $R$ and the ring $S^{-1}R$ is called a {\em maximal left
quotient ring} of the ring $R$ or a {\em maximal left localization
ring} of the ring $R$. The intersection
\begin{equation}\label{llradR}
\gll_R:=\llrad (R) := \bigcap_{S\in \maxDen_l(R)} \ass (S)
\end{equation}
is called the {\em left localization radical } of the ring $R$,
\cite{larglquot}.

 For a ring $R$, there is the canonical exact
sequence 
\begin{equation}\label{llRseq}
0\ra \gll_R \ra R\stackrel{\s }{\ra} \prod_{S\in \maxDen_l(R)}S^{-1}R, \;\; \s := \prod_{S\in \maxDen_l(R)}\, \s_S,
\end{equation}
where $\s_S:R\ra S^{-1}R$, $r\mapsto \frac{r}{1}$. For a left Artinian ring $R$, $|\maxDen_l(R)|<\infty$ and $\gll_R^2= \gll_R$, \cite{Bav-LocArtRing}.


{\bf The maximal elements of $\Ass_l(R)$}.  Let $\maxAss_l(R)$ be
the set of maximal elements of the poset $(\Ass_l(R), \subseteq )$
and

\begin{equation}\label{mADen}
\assmaxDen_l(R) := \{ \ass (S) \, | \, S\in \maxDen_l(R) \}.
\end{equation}
These two sets are equal (Proposition \ref{b27Nov12}), a proof is
based on Lemma \ref{1a27Nov12}. For a non-empty subset $X$ of $R$, let ${\rm r.ass} (X):=\{ r\in R\, | \, rx=0$ for some $x\in X\}$.

\begin{lemma}\label{1a27Nov12}
 \cite{larglquot}
Let $S\in \Den_l(R, \ga )$ and $T\in \Den_l(R, \gb )$ be such that $ \ga \subseteq \gb$. Let $ST$ be the multiplicative semigroup generated by $S$ and $T$ in $(R,\cdot )$.  Then
\begin{enumerate}
\item ${\rm r.ass} (ST)\subseteq \gb$.
\item $ST \in \Den_l(R, \gc )$ and $\gb \subseteq \gc$.
\end{enumerate}
\end{lemma}


\begin{proposition}\label{b27Nov12}
\cite{larglquot} $\; \maxAss_l(R)= \assmaxDen_l(R)\neq \emptyset$. In particular, the ideals of this set are incomparable (i.e. neither $\ga\nsubseteq \gb$ nor $\ga\nsupseteq \gb$).
\end{proposition}

{\bf The localization maximal rings}. The set ($\Loc_l(R), \ra )$ is a partially ordered set (poset) where $A\ra B$ if there is an $R$-homomorphism $A\ra B$. There is no oriented loops in the poset $\Loc_l(R)$ apart from the $R$-isomorphism $A\ra A$.

Let $\maxLoc_l(R)$ be the set of maximal elements of the poset
$(\Loc_l(R), \ra )$. Then (see \cite{larglquot}),
\begin{equation}\label{mADen1}
\maxLoc_l(R) = \{ S^{-1}R \, | \, S\in \maxDen_l(R) \}= \{ Q_l(R/
\ga ) \, | \, \ga \in \assmaxDen_l(R)\}.
\end{equation}


 {\it Definition}, \cite{larglquot}. A ring $A$ is
called a {\em left localization maximal ring} if $A= Q_l(A)$ and
$\Ass_l(A) = \{ 0\}$. A ring $A$ is called a {\em right
localization maximal ring} if $A= Q_r(A)$ and $\Ass_r(A) = \{
0\}$. A ring $A$ which is a left and right localization maximal
ring is called a {\em (left and right) localization maximal ring}
(i.e. $Q_l(A) =A=Q_r(A)$ and $\Ass_l(A) =\Ass_r(A) = \{ 0\}$).

$\noindent $

The next theorem is a criterion of  when a left quotient ring of a
ring is a maximal left quotient ring of the ring.

\begin{theorem}\label{21Nov10}
\cite{larglquot} Let  a ring $A$ be a left localization of a ring
$R$, i.e. $A\in \Loc_l(R, \ga )$ for some $\ga \in \Ass_l( R)$.
Then $A\in \maxLoc_l( R)$ iff $Q_l( A) = A$ and  $\Ass_l(A) = \{
0\}$, i.e. $A$ is a left localization maximal ring (clearly, $\ga \in \assmaxDen_l(R)$).
\end{theorem}


Theorem \ref{21Nov10} shows that the left localization maximal
rings are precisely the localizations of all the rings at their
maximal left denominators sets.


{\it Example}. Let $A$ be a simple ring. Then $Q_l(A)$ is a left
localization maximal  ring and $Q_r(A)$ is a right localization
maximal ring (by Theorem \ref{4Jul10}.(4) and Theorem \ref{21Nov10}).


{\it Example}. A division ring is a (left and right) localization
maximal ring. More generally, a simple Artinian ring (i.e. the
matrix algebra over a division ring) is a (left and right)
localization maximal ring.


{\bf Left (non-)localizable elements of a ring}.

\begin{lemma}\label{b11Dec12}
\cite{Bav-LocArtRing} Let $S\in \Den_l(R, \ga )$ (respectively, $S\in \Den (R, \ga )$),
$\s : R\ra S^{-1}R$, $r\mapsto \frac{r}{1}$, and $G:=(S^{-1}R)^*$
be the group of units of the ring $S^{-1}R$. Then $S':=\s^{-1}
(G)\in \Den_l(R, \ga )$ (respectively, $S':=\s^{-1} (G)\in \Den
(R, \ga )$).
\end{lemma}


{\it Definition}, \cite{Crit-S-Simp-lQuot}. An element $r$ of a ring $R$ is called a {\em left localizable element} if there  exists a left denominator set $S$  of
 $R$ such that $r\in S$ (and so the element $\frac{r}{1}\neq 0$ is invertible in the ring
 $S^{-1}R$), equivalently,  if there  exists a left denominator set $T$  of
 $R$ such that the element $\frac{r}{1}$ is invertible in the ring
 $T^{-1}R$ (Lemma \ref{b11Dec12}). The set of left localizable elements is denoted $\CL_l(R)$.

$\noindent $

Clearly, 
\begin{equation}\label{CLeU}
\CL_l(R)=\bigcup_{S\in \maxDen_l(R)} S.
\end{equation}





{\bf The maximal left quotient rings of a finite direct product of rings}.
\begin{theorem}\label{c26Dec12}
\cite{Crit-S-Simp-lQuot}  Let $R=\prod_{i=1}^n R_i$ be a direct product of rings $R_i$. Then
for each $i=1, \ldots , n$, the map
\begin{equation}\label{1aab1}
\maxDen_l(R_i) \ra \maxDen_l(R), \;\; S_i\mapsto R_1\times\cdots \times S_i\times\cdots \times R_n,
\end{equation}
is an injection. Moreover, $\maxDen_l(R)=\coprod_{i=1}^n \maxDen_l(R_i)$ in the sense of (\ref{1aab1}), i.e.
$$ \maxDen_l(R)=\{ S_i\, | \, S_i\in \maxDen_l(R_i), \; i=1, \ldots , n\},$$
$S_i^{-1}R\simeq S_i^{-1}R_i$, $\ass_R(S_i)= R_1\times \cdots \times \ass_{R_i}(S_i)\times\cdots \times R_n$. The core of the left denominator set $S_i$ in $R$ coincides with the core $S_{i,c}$ of the left denominator set $S_i$ in $R_i$, i.e.
$$(R_1\times\cdots \times 0\times S_i\times 0\times\cdots \times R_n)_c=0\times\cdots \times 0\times S_{i,c}\times 0\times\cdots \times 0.$$
\end{theorem}


{\bf Properties of the maximal left quotient rings of a ring}.
The next theorem describes various properties of the maximal left
quotient rings of a ring, in particular, their groups of units and
their largest left quotient rings.

\begin{theorem}\label{15Nov10}
\cite{larglquot} Let $S\in \maxDen_l(R)$, $A= S^{-1}R$, $A^*$ be
the group of units of the ring $A$; $\ga := \ass (S)$, $\pi_\ga
:R\ra R/ \ga $, $ a\mapsto a+\ga$, and $\s_\ga : R\ra A$, $
r\mapsto \frac{r}{1}$. Then
\begin{enumerate}
\item $S=S_\ga (R)$, $S= \pi_\ga^{-1} (S_0(R/\ga ))$, $ \pi_\ga
(S) = S_0(R/ \ga )$ and $A= S_0( R/\ga )^{-1} R/ \ga = Q_l(R/ \ga
)$. \item  $S_0(A) = A^*$ and $S_0(A) \cap (R/ \ga )= S_0( R/ \ga
)$. \item $S= \s_\ga^{-1}(A^*)$. \item $A^* = \langle \pi_\ga (S)
, \pi_\ga (S)^{-1} \rangle$, i.e. the group of units of the ring
$A$ is generated by the sets $\pi_\ga (S)$ and $\pi_\ga^{-1}(S):=
\{ \pi_\ga (s)^{-1} \, | \, s\in S\}$. \item $A^* = \{ \pi_\ga
(s)^{-1}\pi_\ga ( t) \, |\, s, t\in S\}$. \item $Q_l(A) = A$ and
$\Ass_l(A) = \{ 0\}$.     In particular, if $T\in \Den_l(A, 0)$
then  $T\subseteq A^*$.
\end{enumerate}
\end{theorem}



{\bf A bijection between $\maxDen_l(R)$ and $\maxDen_l(Q_l(R))$}.
\begin{proposition}\label{A8Dec12}
\cite{Crit-S-Simp-lQuot} Let $R$ be a ring, $S_l$ be the  largest regular left Ore set of the ring $R$, $Q_l:= S_l^{-1}R$ be the largest left quotient ring of the ring $R$, and $\CC$ be the set of regular elements of the ring $R$. Then
\begin{enumerate}
\item $S_l\subseteq S$ for all $S\in \maxDen_l(R)$. In particular,
$\CC\subseteq S$ for all $S\in  \maxDen_l(R)$ provided $\CC$ is a
left Ore set. \item Either $\maxDen_l(R) = \{ \CC \}$ or,
otherwise, $\CC\not\in\maxDen_l(R)$. \item The map $$
\maxDen_l(R)\ra \maxDen_l(Q_l), \;\; S\mapsto SQ_l^*=\{ c^{-1}s\,
| \, c\in S_l, s\in S\},
$$ is a bijection with the inverse $\CT \mapsto \s^{-1} (\CT )$
where $\s : R\ra Q_l$, $r\mapsto \frac{r}{1}$, and $SQ_l^*$ is the
sub-semigroup of $(Q_l, \cdot )$ generated by the set  $S$ and the
group $Q_l^*$ of units of the ring $Q_l$, and $S^{-1}R= (SQ_l^*)^{-1}Q_l$.
    \item  If $\CC$ is a left Ore set then the map $$ \maxDen_l(R)\ra \maxDen_l(Q), \;\; S\mapsto SQ^*=\{ c^{-1}s\,
| \, c\in \CC, s\in S\}, $$ is a bijection with the inverse $\CT
\mapsto \s^{-1} (\CT )$ where $\s : R\ra Q$, $r\mapsto
\frac{r}{1}$, and $SQ^*$ is the sub-semigroup of $(Q, \cdot )$
generated by the set  $S$ and the group $Q^*$ of units of the ring
$Q$, and $S^{-1}R= (SQ^*)^{-1}Q$.
\end{enumerate}
\end{proposition}










\section{Weakly left localizable rings and their
characterizations}\label{WLLCH}

In this section, a new class of rings is introduced, the {\em
weakly left localizable rings}, and their characterizations are given
(Theorem \ref{26Mar14} and Theorem \ref{28Mar14}) under certain natural conditions.

$\noindent $

{\it Definition}. A ring $R$ is called a {\em weakly left
localizable ring} (a WLL ring, for short) if  every {\em
non-nilpotent} element is left localizable, i.e. the ring
\begin{equation}\label{R=LN}
R=\CL_l(R)\coprod \Nil (R)
\end{equation}
is a disjoint union of the set $\CL_l:=\CL_l(R)$ of left
localizable elements and the set $N:=\Nil (R)$ of nilpotent
elements of the ring $R$.

$\noindent $

Every left localizable ring is a weakly left localizable ring but not vice versa. An ideal of a ring is called a {\em nil ideal} if every element of
it is nilpotent. The sum and the product of two nil ideals is a
nil ideal. The sum $\CN = \CN_R$ of all nil ideals of a ring $R$
is called the {\em nil radical} of the ring $R$. The nil radical
$\CN$ is the largest nil ideal of the ring $R$. Clearly,
$$ \gn_R\subseteq \CN_R\subseteq \rad (R) \;\; {\rm and}\;\; \CN_R\subseteq
\Nil (R)$$
where $\gn_R$ is the {\em prime radical} of the ring $R$.
The next lemma provides many examples of  weakly left localizable rings.

\begin{lemma}\label{a8Feb13}
\begin{enumerate}
\item Every left localizable ring is weakly left localizable but not vice versa.
\item Let $R$ be a ring with $R = R^*\coprod \Nil (R)$. Then $R$ is a weakly left localizable ring.
\item Every left Artinian ring $R$ such that $R/ \gn$ is a division ring is a weakly left localizable ring (where $\gn$ is the prime radical of $R$).
\item  Every local ring $R$ such that $\rad (R) = \CN_R$ is a weakly left localizable ring.
\end{enumerate}
\end{lemma}

{\it Proof}. 1 and 2. Obvious.

3. Since $R=R^*\coprod \Nil (R)$, $R$ is a weakly left localizable ring, by statement 2.

4. Clearly, $\rad (R) = \CN_R = \Nil (R)$. Then $R= R^*\cup \rad (R)= R^* \cup \Nil (R)$, i.e. $R$ is a weakly left localizable ring. $\Box $

$\noindent $

\begin{theorem}\label{9Feb13}
Let $R=\prod_{i=1}^n R_i$ be a direct product of rings $R_i$. The ring $R$ is a weakly left localizable ring iff the rings $R_i$ are so.
\end{theorem}

{\it Proof}. The proof is an easy corollary of Theorem \ref{c26Dec12} that states that $\maxDen_l(\prod_{i=1}^n R_i)= \coprod_{i=1}^n \maxDen_l(R_i)$.

$(\Rightarrow )$ Suppose that the ring $R$ is a weakly  left localizable ring. We have to show that the rings $R_i$ are so. Each  ring $R_i$ is a subring of $R$. Let $r_i\in R_i$ be a non-nilpotent  element. The ring $R_i$ is a weakly  left localizable ring. So, $r_i\in S_i$ for some $S_i\in \maxDen_l(R)$. Then $S_i\in \maxDen_l(R_i)$, by Theorem \ref{c26Dec12}.

$(\Leftarrow )$ Suppose that the rings $R_i$ are weakly  left localizable rings. Let $r= (r_1, \ldots , r_n)\in R$ be a non-nilpotent  element. Then $0\neq r_i\in R_i$ is  a non-nilpotent  element for some $i$. The ring $R_i$ is a weakly left localizable ring, and so $r_i\in S_i$ for some $S_i\in \maxDen_l(R_i)$. By Theorem \ref{c26Dec12}, $S_i\in \maxDen_l(R)$. Therefore, $R$ is a weakly left localizable ring. $\Box $

$\noindent $

{\it Definition}, \cite{Crit-S-Simp-lQuot}. For an arbitrary ring $R$, the intersection $$\CC_l(R):=\bigcap_{S\in \maxDen_l(R)}S$$
is called the set of {\em completely left localizable elements} of $R$ and an element of the set $\CC_l(R)$ is called a {\em completely left  localizable element}.

\begin{lemma}\label{a28Mar14}
 Let $R$ be a ring, $\gll := \gll_R$, $\pi': R\ra R':= R/\gll$, $r\mapsto \br := r+\gll$; $\CC_l(R)$ and $\CC_l(R')$ be the sets of completely left localizable elements of the rings $R$ and $R'$ respectively. Then
\begin{enumerate}
\item  the map $\phi : \maxDen_l(R)\ra \maxDen_l(R')$, $S\mapsto \pi'(S)$, is an injection.
\item $\pi'^{-1}(\CC_l(R'))\subseteq \CC_l(R)$.
\item If the map $\phi$ is a surjection then $\pi'^{-1} (\CC_l(R'))= \CC_l(R)$ and $\pi' (\CC_l(R))= \CC_l(R')$.
\end{enumerate}
\end{lemma}

{\it Proof}. 1. The inclusion  $\gll\subseteq \ass (S)$ implies that  $\pi'(S)\in \Den_l(R' , \ass (S)/\gll )$
 and $S^{-1}R \simeq \pi'(S)^{-1}R'$ is a left localization maximal ring. Hence, $\pi'(S)\in \maxDen_l(R')$ since $S+\ass (S) \subseteq S$ and   $\gll\subseteq \ass (S)$. Since $S+\ass (S) \subseteq S$, we have the inclusion $S+\gll \subseteq S$ for all $S\in \maxDen_l(R)$, i.e. the map $S\mapsto \pi'(S)$ is an injection.

2. Statement 2 follows from statement 1 and the fact that $S+\gll \subseteq S$ for all $S\in \maxDen_l(R)$. In more detail,
$$ \pi'^{-1}(\CC_l(R'))= \bigcap_{S'\in \maxDen_l(R')}\pi'^{-1}(S')\subseteq
 \bigcap_{S\in \maxDen_l(R)}\pi'^{-1}(\pi'(S))=\bigcap_{S\in \maxDen_l(R)}S=\CC_l(R).$$

3. The equality $\pi'^{-1}(\CC_l(R'))= \CC_l(R)$ is obvious (use the proof of statement 2 where the inclusion is replaced by the equality).
\begin{eqnarray*}
\pi'(\CC_l(R))&=& \pi' (\bigcap_{S\in \maxDen_l(R)}S)=
 \bigcap_{S\in \maxDen_l(R)}\pi'(S) \;\;\; ({\rm since}\;\; S+\gll \subseteq S )\\
 &=&\bigcap_{S'\in \maxDen_l(R')}S'=\CC_l(R'),
\end{eqnarray*}
since $\phi$ is a surjection.  $\Box $

The next theorem shows that the set of maximal denominator sets behaves nicely under localizations at regular left denominator sets, it is used in the proof of Proposition \ref{a14Dec12}.


\begin{theorem}\label{C3Dec12}
\cite{left-localizable-rings} Let $R$ be a ring, $T\in \Den_l(R, 0)$ and $\s : R\ra T^{-1}R$,
$r\mapsto \frac{r}{1}$. Then
\begin{enumerate}
\item $T\subseteq S$ for all $S\in \maxDen_l(R)$.  \item  The map
$$ \maxDen_l(R)\ra \maxDen_l(T^{-1}R), \;\; S\mapsto
\widetilde{S},$$ is a bijection with the inverse $\CT \mapsto
\s^{-1} (\CT )$ where $\widetilde{S}$ is the multiplicative monoid
generated in the ring $T^{-1}R$ by $\s (S)$ and $\s (T)^{-1}:=\{
t^{-1} \, | \, t\in T\}$, and $S^{-1}R\simeq \widetilde{S}^{-1}(T^{-1}R)$.
\end{enumerate}
\end{theorem}

\begin{proposition}\label{a14Dec12}
Let $R$ be a ring. Then
\begin{enumerate}
\item $\gll_{R/ \gll_R}=0$.
\item $\CL_l(R)+\gll_R\subseteq \CL_l(R)$.
\item $\CL_l(R/ \gll_R)\supseteq  \CL_l(R)+ \gll_R$.
\item Let $S\in \Den_l(R, 0)$. Then $\gll_{S^{-1}R}=S^{-1}\gll_R$. In particular, $\gll_{Q_l(R)}= S_l(R)^{-1}\gll_R$ and $\gll_{Q_{l, cl}(R)}= \CC_R^{-1}\gll_R$.
\end{enumerate}
\end{proposition}

{\it Proof}. 1. We keep the notation of Lemma \ref{a28Mar14}. Statement 1 follows from Lemma \ref{a28Mar14}.(1),
\begin{eqnarray*}
\gll_{R/ \gll_R} &\subseteq &\bigcap_{S\in \maxDen_l(R)}\ass (\pi'(S))= \bigcap_{S\in \maxDen_l(R)}\ass (S)/\gll_R \\
 &=& (\bigcap_{S\in \maxDen_l(R)}\ass (S))/\gll_R =\gll_R/\gll_R=0.
\end{eqnarray*}
2 and 3.  Statements 2 and 3 follow from  (\ref{llRseq}).

4. By Theorem \ref{C3Dec12}.(2), $S^{-1}R\simeq \widetilde{S}^{-1}(T^{-1}R)$ for all $S\in \maxDen_l(R)$. So, there is a natural ring isomorphism $\tau : Q:=\prod_{S\in \maxDen_l(R)}S^{-1}R\ra Q':=\prod_{S\in \maxDen_l(R)} \widetilde{S}^{-1}(T^{-1}R)$. As a result we have the commutative diagram

\[\renewcommand{\arraystretch}{1.5}
\begin{array}{ccccc}
R &\xrightarrow{\s }&Q \\
\Big\downarrow{\tau'} &&\Big\downarrow{\tau}\\
S^{-1}R &\xrightarrow{\s'  }& Q'
\end{array}\]
where $\tau'$ is a monomorphism given by the rule $\tau' (r) = \frac{r}{1}$ and $\s$, $\s'$ are defined in (\ref{llRseq}). Let $s^{-1}r\in S^{-1}R$ where $s\in S$ and $r\in R$. Then $s^{-1}r\in \gll_{S^{-1}R}$, i.e. $0=\s'(s^{-1}r)$, iff $0=\s'\tau'(r)= \tau \s (r)$ iff $\s (r)=0$ iff $r\in \gll_R$. Therefore, $\gll_{S^{-1}R}=S^{-1}\gll_R$.  $\Box $


{\bf Weakly left localizability criterion}. The next proposition is a weakly left localizability criterion.
\begin{proposition}\label{c13Dec12}
A ring $R$ is a weakly left localizable ring  iff the factor ring $R/ \gll_R$ is a
weakly left localizable ring,  the left localization radical $\gll_R$ of $R$  is a nil ideal  and $\pi' (\CL_l(R))= \CL_l(R/ \gll_R)$ where $\pi' : R\ra R/ \gll_R$, $r\mapsto \br := r+\gll_R$.
\end{proposition}

{\it Proof}. $(\Rightarrow )$ The ring $R$ is a weakly left localizable ring, hence $\gll_R\subseteq \Nil (R)$, by (\ref{R=LN}), i.e. $\gll_R$ is a nil ideal. Then, by Proposition \ref{a14Dec12}.(2,3) and (\ref{R=LN}),
\begin{equation}\label{R=LN1}
\pi'(R)= \pi'(\CL_l(R)) \coprod \pi'(\Nil (R)).
\end{equation}
 An element $r\in R$ is nilpotent iff the element $r+\gll_R\in R/ \gll_R$ is nilpotent (since $\gll_R$ is a nil ideal). Therefore, $\pi' (\Nil (R))= \Nil (R/ \gll_R))$.  Since $\pi' (\CL_l(R))\subseteq \CL_l(R/ \gll_R)$, we must have $\pi' (\CL_l (R))=\CL_l (R/ \gll_R))$, by (\ref{R=LN1}). Now, $R/ \gll_R=\CL_l(R/ \gll_R)\coprod\Nil (R/ \gll_R)$, by (\ref{R=LN1}). So, $R/ \gll_R$ is a weakly left localizable ring.

$(\Leftarrow )$ Let $r\in R$ and $\gll := \gll_R$. Then $r$ is a
nilpotent element in $R$ iff $r+ \gll$ is is a nilpotent element
in the factor ring $R/ \gll$ since $\gll$ is a nil ideal. Suppose that  $r$
is not a nilpotent element of the ring $R$. Then necessarily $\br\in \CL_l(R/ \gll_R)$ and $\br= \br_1$ for some $r_1\in \CL_l(R)$ (since $\pi' (\CL_l(R))= \CL_l(R/ \gll_R)$). By Proposition \ref{a14Dec12}.(2), $r\in r_1+\gll_R\in \CL_l(R)$. Therefore, $R$ is a weakly left localizable ring.  $\Box $



\begin{lemma}\label{d13Dec12}
Let $R$ be a ring. If $yx=1$ for some elements $x,y\in R$ then
neither $x$ nor $y$ is a nilpotent element.
\end{lemma}

{\it Proof}. If $x^n=0$ (respectively, $y^n=0$) for some $n\geq 1$ then
$1=y^nx^n=0$ (respectively, $1=y^nx^n=0$), a contradiction. $\Box
$


\begin{lemma}\label{a26Mar14}
Let $R$ be a ring. Then $R$ is a left localization maximal and weakly left localizable ring (i.e. $R=R^*\coprod \Nil (R)$) iff $\Nil (R) = \CN_R$ and $R=R^*\coprod \CN_R$, i.e. $R$ is a local ring  with $\rad (R) = \CN_R$.
\end{lemma}

{\it Proof}. $(\Rightarrow )$ Let $N=\Nil (R)$ and $\CN = \CN_R$. It suffices to show that  $N = \CN$ since $R= R^*\coprod N = R^*\coprod \CN$.  To prove that $N=\CN$ it suffices to show that $N$ is an ideal of $R$ (since $N\supseteq \CN$ and $\CN$ is the largest nil ideal of $R$).

(i) $R^*NR^*\subseteq N$: If this inclusion does not hold, i.e. $w:= unv\in R^*$ for some elements $u,v\in R^*$ and $n\in N$, then $nvw^{-1} u=1$, a contradiction, by Lemma \ref{d13Dec12}.

(ii) $NN\subseteq N$: If this inclusion does not hold, i.e. $u:= n_1n_2\in R^*$ for some $n_1, n_2\in N$ then $u^{-1}n_1\cdot n_2=1$, a contradiction, by Lemma \ref{d13Dec12}.

(iii) $RNR\subseteq N$: By (i), (ii) and $R= R^*\coprod N$.

(iv) $N+N\subseteq N$: If this inclusion does not hold, i.e. $u:= n_1+n_2\in R^*$ for some elements $n_1, n_2\in N$. Then, by (i), $N\ni u^{-1} n_1= 1-u^{-1} n_2\in 1-N \subseteq R^*$, a contradiction.

(v) $N$ {\em is an ideal}, by (iii) and (iv).

$(\Leftarrow )$ This implication is obvious.  $\Box $


{\bf The left quotient rings of a class of weakly left localizable rings}.
 The next theorem (Theorem \ref{26Mar14}) is  a criterion (via weakly left localizable rings) for a ring to have the classical left quotient ring to be a direct product $\prod_{i=1}^n R_i$ of local rings with $\rad (R_i) = \CN_{R_i}$ for $i=1, \ldots , n$. It is also a description of a certain class of weakly left localizable rings.

\begin{theorem}\label{26Mar14}
Let $R$ be a ring. The following statements are equivalent.
\begin{enumerate}
\item The ring $R$ is a weakly left localizable ring such that
\begin{enumerate}
\item $\gll_R=0$,
\item $|\maxDen_l(R)|<\infty$,
\item for every $S\in \maxDen_l(R)$, $S^{-1}R$ is a weakly left localizable ring, and
    \item  for all $S, T\in \maxDen_l(R)$ such that $S\neq T$, $\ass (S)$ is not a nil ideal modulo $\ass (T)$.
\end{enumerate}
\item $Q_{l, cl}(R)\simeq \prod_{i=1}^n R_i$ where $R_i$ are local rings with $\rad (R_i) = \CN_{R_i}$.
\item $Q_l(R)\simeq \prod_{i=1}^n R_i$ where $R_i$ are local rings with $\rad (R_i) = \CN_{R_i}$.
\end{enumerate}
\end{theorem}

{\it Proof}. $(1 \Rightarrow 2)$ By the condition (b), $\maxDen_l(R)=\{ S_1, \ldots , S_n\}$. Let $\ga_i := \ass (S_i)$. For each $i=1, \ldots , n$, the ring $R_i:= S_i^{-1}R$ is a left localization maximal ring (Theorem \ref{21Nov10}) and a weakly left localization maximal  ring, by the condition (c). By Lemma \ref{a26Mar14}, we have the statement (i) below.

(i) $R_i=R_i^*\coprod \CN_i$ {\em  where} $\CN_i :=\CN_{R_i}=\Nil (R_i)$:

If $n=1$ then we are done, i.e. the implication $(1\Rightarrow 2)$ is true. In more detail, if $n=1$ then $\ass (S_1) = \gll =0$, and so $S_1\subseteq \CC_R$. The left $R$-module $R$ is an essential submodule of $R_1:= S_1^{-1}R$. Hence, for every element $c\in \CC_R$, the map $\cdot c : R_1\ra R_1$, $a\mapsto ac$, is an injection. By (i), $c\in R_1^*$, i.e. $ \CC_R\subseteq R_i^*$. This implies that $\CC_R\in \Den_l(R, 0)$ with $Q_{l,cl} (R)=R_1$.

Therefore, we can assume that $n\geq 2$. By the statement (a), the map
$$\s : R\ra \prod_{i=1}^n R_i, \;\; r\mapsto (r_1, \ldots , r_n), \;\; r_i:= \frac{r}{1}\in R_i,$$ is a monomorphism and we can identify the ring $R$ with its image $\s (R)$.

(ii) {\em For} $i=1, \ldots , n$, $\CC_i':= S_i \cap \cap_{j\neq i}\ga_j\neq \emptyset$: We may assume that $i=1$. Fix $r\in S_1$. Then $r=(r_1, \ldots , r_n) \in \prod_{i=1}^n R_i$ with $r_1\in R_1^*$. Replacing $r$ by $r^t$ for some $t\geq 1$ we may assume that $r_j\in R_j^*\cup \{ 0\}$ (by (i)) for all $j=2,\ldots , n$. Fix an element $r\in S_1$ with the least number of non-zero coordinates. Up to order, we may assume that $r=(r_1, \ldots , r_s, 0,\ldots , 0)$ where $r_i\in R_i^*$ for $i=1,  \ldots , s$. We claim that $s=1$. Suppose that $s>1$, we seek a contradiction. Notice that
$$\ga_sr= (\ga_sr_1, \ldots , \ga_sr_{s-1}, 0,\ldots , 0).$$ We claim that $\ga_sr_1\cap R_1^*\neq \emptyset$. Suppose that  $\ga_sr_1\cap R_1^*=\emptyset$, then $\ga_sr_1\subseteq \CN_1$ (by (i)). Hence, $\frac{\ga_s}{1} \subseteq \CN_1r^{-1}= \CN_1$, and so the ideal $\ga_s$ is a nil ideal modulo $\ga_1$, this contradicts to the condition (d). Therefore, $\frac{a_s}{1}r_1\in R_1^*$ for some element $a_s\in \ga_s$, and so $a_sr_1\in \s_1^{-1}(R_1^*)=S_1$ where $\s_1: R\ra R_1$, $ x\mapsto \frac{x}{1}$ (by Theorem \ref{15Nov10}.(3)). This contradicts to the minimality of $s$. So, $s=1$ and therefore the element $r= (r_1, 0 ,\ldots , 0)\in \CC_1'$ since $r\in R = (\cup_{i=1}^n S_i) \coprod\Nil (R)$ (as $R$ is a weakly left localizable ring) and $r\not\in (\cup_{i=2}^n S_i) \coprod\Nil (R)$, by the choice of the element $r= (r_1, 0 , \ldots , 0)$.

(iii) $\CC_i'\CC_j'=0$ {\em for all} $i\neq j$:  $\CC_i'\CC_j'=\cap_{k=1}^n \ga_k = \gll_R =0$, by the condition (a).

(iv) $S_i\CC_i'S_i\subseteq \CC_i'$ {\em for} $i=1, \ldots , n$: Obvious, by (i).

(v) $\ass (\CC_i') := \{ r\in R\, | \, c_i'r=0$ {\em for some} $c_i'\in\CC_i'\}=\ga_i$: The inclusion $ \CC_i'\subseteq S_i$ implies the inclusion $\ass (\CC_i') \subseteq \ga_i$. The reverse inclusion follows from the inclusion $\CC_i'S_i\subseteq \CC_i'$ (see (iv)).

(vi) $\CC_i'\in \Den_l(R, \ga_i)$ {\em and}  $\CC_i'^{-1}R\simeq S_i^{-1}R$: Clearly, the set $\CC_i'$ is a multiplicative set not necessarily containing 1. Fix an element $c_i\in \CC_i'$. Then  each element $s^{-1}r\in R_i$, where $s\in S_i$ and $r\in R$, can be written as a left fraction
\begin{equation}\label{s1rC}
s^{-1}r= s^{-1} c_i^{-1}c_i r = (c_is)^{-1} c_ir\;\; {\rm where} \;\; c_is\in \CC_i'\;\; {\rm and}\;\; c_ir \in \cap_{j\neq i} \ga _j,
\end{equation}
by (iv). The the statement (vi) follows bearing in mind (v).

(vii) $\CC':= \CC_1'+\cdots +\CC_n'\in \Den_l(R, 0)$ {\em and} $\CC'^{-1} R \simeq \prod_{i=1}^n R_i$: Since $\s$ is a monomorphism, $\CC'\subseteq \CC_R$. By (\ref{s1rC}), any element $\alpha$ of the ring $\prod_{i=1}^n R_i$ can be written as $\alpha = (c_1^{-1}r_1. \ldots , c_n^{-1}r_n)$ where $c_i\in \CC_i'$ and $r_i\in \cap_{j\neq i}\ga_j$. Hence, $c_ic_j=0$ for all $i\neq j$, by (iii). Then $c:=c_1+\cdots + c_n\in \CC'$, $r:= r_1+\cdots + r_n\in R$ and
\begin{eqnarray*}
c\alpha &= & (cc_1^{-1}r_1, \ldots , cc_n^{-1}r_n) =  (cc_1c_1^{-2}r_1, \ldots , cc_nc_n^{-2}r_n)\\
 &=& (c_1^2c_1^{-2}r_1, \ldots , c_n^2c_n^{-2}r_n)= (r_1, \ldots , r_n)= r_1+\cdots + r_n = r
\end{eqnarray*}
since $r_i\in \cap_{j\neq i} \ga_j$. So, $\alpha = c^{-1} r$, and the statement (vii) follows.

(viii) $Q_{l, cl} (R) \simeq \CC'^{-1}R$: Let $Q':= \CC'^{-1}R$. Recall that $Q'=\prod_{i=1}^n R_i$, $R_i = R_i^*\coprod S_i^{-1}\CN$ and $\CC'\subseteq \CC_R$. The ring $R$ is an essential left $R$-submodule of $Q'$. Therefore, for each $c\in \CC_R$, the $R$-module homomorphism $\cdot c:Q'\ra Q'$, $q' \mapsto q'c$, is a monomorphism. Then $c\in \prod_{i=1}^nR_i^*= Q'^*$. So, $\CC' \subseteq \CC_R\subseteq Q'^*$. Since $Q' = \CC'^{-1}R$ we must have $\CC_R\in \Den_l(R)$ and $Q'= \CC_R^{-1}R = Q_{l, cl} (R)$.

$(2\Rightarrow 3)$ Trivial.

$(3\Rightarrow 1)$ Notice that $R_i= R_i^*\coprod\CN_{R_i}$, and so $R_i$ is a weakly left localizable ring and $\maxDen_l(R_i) = \{ R_i^*\}$. By Theorem \ref{c26Dec12},
$$ \maxDen_l(Q_l(R)) = \maxDen_l(\prod_{i=1}^n R_i) = \{ S_i'\, | \, i=1, \ldots , n\}$$ where $S_i':= R_1\times\cdots \times R_{i-1}\times R_i^*\times R_{i+1}\times \cdots \times R_n$ and $\ass (S_i')= R_1\times\cdots \times R_{i-1}\times \{ 0\} \times R_{i+1}\times \cdots \times R_n$. By Proposition \ref{A8Dec12},
$ \maxDen_l(R)=\{ S_i\, | \, i=1, \ldots , n\}$ is a finite set where $S_i=\s'^{-1}(S_i')$ and $\s' : R\ra Q_l(R)$, $r\mapsto\frac{r}{1}$; and $ S_i^{-1} R\simeq S_i'^{-1}Q_l(R)=R_i$ is a weakly left localizable ring. It follows from the commutative diagram
\[\renewcommand{\arraystretch}{1.5}
\begin{array}{ccccc}
R &\xrightarrow{\s }&\prod_{i=1}^nS_i^{-1}R \\
&\searrow{\s'}&\Big\downarrow\rlap{{\rm id}}\\
&& Q_l(R)=\prod_{i=1}^nR_i
\end{array}\]
 that $\gll_R=0$ since $\s'$ is a monomorphism where id is the identity map. If $r\not\in \Nil (R)$ then necessarily $r\in S_i'$ for some $i$, hence $r\in \s'^{-1}(S_i')= S_i$, i.e. the ring $R$ is a weakly left localizable ring.

To prove that the statement (d) holds we may assume that $n>1$. It suffices to show that the ideal $\ga_2$ is not a nil ideal modulo $\ga_1$. The element $e_1:= (1, 0 , \ldots , 0)\in R_1^*\cap \ass (S_2')$ can be written as $c^{-1}r$ for some elements $c\in  \CC_R$ and $r\in R$. Then the element $r= ce_1\in R\cap R_1^*\cap \ass (S_2')= S_1\cap \ga_2$ is obviously not a nilpotent element modulo $\ga_1$.  $\Box $

$\noindent $

The next corollary shows that the rings that satisfy the conditions of Theorem \ref{26Mar14} have interesting properties. In particular, for each such a ring $R$, the set of completely left localizable elements coincide with the set of regular elements, $\CC_R= \CC_l(R)$,  and $\Nil (R) = \CN_R$. Notice that for an arbitrary ring $R$ the inclusions hold, $\CC_R\subseteq \CC_l(R)$ (by Proposition \ref{A8Dec12}.(1), provided $\CC_R\in \Ore_l(R)$) and $\Nil (R) \supseteq  \CN_R$, that are, in general, are not equalities.
\begin{corollary}\label{b26Mar14}
We keep the notation of Theorem \ref{26Mar14}. Suppose that a ring $R$ satisfies one of the equivalent conditions 1--3 of Theorem \ref{26Mar14}. Then
\begin{enumerate}
\item $\maxDen_l(R)=\{ S_1, \ldots , S_n\}$ where $S_i=\{ r\in R\, | \, \frac{r}{1}\in R_i^*\}$.
\item $\CC_R = \CC_l(R)$, i.e. $\CC_R=\bigcap_{S\in \maxDen_l(R)}S$.
\item $\Nil (R) = \CN_R$.
\item $Q:= Q_{l,cl} (R) = Q_l(R)$ is a weakly left localizable ring with $\Nil (Q) = \CN_Q= \rad (Q)$.
    \item $\CC_R^{-1}\CN_R =\CN_Q
     = \rad (Q)$ and $\CN_R=R\cap \CN_Q$.
    \item $\CC_R^{-1} \CL_l(R) = \CL_l(Q)$ and $\CL_l(R)=R\cap \CL_l(Q)$.
\end{enumerate}
\end{corollary}

{\it Proof}. 1. Statement 1 has already been established in the proof of Theorem \ref{26Mar14}.

2. By  Proposition \ref{A8Dec12}.(1), $\CC_R\subseteq \CC_l(R)=\bigcap_{S\in \maxDen_l(R)}S$. The inverse inclusion follows from the fact that $\CC_l(R)\subseteq \prod_{i=1}^n R_i^* = Q_{l,cl}(R)^*$.

3. Let $r\in R$. It follows from the inclusion $R\subseteq Q_{l,cl}(R) = \prod_{i=1}^n R_i$ that $r\in N$ iff $r\in \prod_{i=1}^n \CN_{R_i}= \CN_{Q_{l,cl}(R)}$. It follows that $N$ is an ideal, hence $N=\CN_R$.

4. The equality $Q_{l,cl}(R) = Q_l(R)$ has already been established. By Theorem \ref{26Mar14}.(2), $\Nil (Q) = \CN_Q= \rad (Q)$.

5 and 6. By statement 4, $Q= \CL_l(Q)\coprod \CN_Q$. Since $\CN_R\subseteq \CN_Q$, we have $\CC_R^{-1} \CN_R\subseteq \CN_Q$. Since $\CC_R\subseteq S$ for all maximal left denominator sets $S\in \maxDen_l(R)$, we have the equality $\CC_R^{-1}\CL_l(R):= \{ c^{-1} s\, | \, c\in \CC_R, s\in \CL_l(R)\}\subseteq \CL_l(Q)$. Then the first equalities of  the statements 5 and 6 follow. Then, $\CN_R\subseteq R\cap \CC_R^{-1}\CN_R=R\cap \CN_Q\subseteq \CN_R$, i.e. $\CN_R=R\cap \CN_Q$. Finally, the equality  $\CL_l(R)=R\cap \CL_l(Q)$ is obvious.  $\Box $

$\noindent $

The next proposition describes the set $\CC_l(R)$ (statement 2) and gives sufficient conditions  for $\CC_l(R)$ to be a left denominator set in $R$.
\begin{proposition}\label{b28Mar14}
We keep the notation of Lemma \ref{a28Mar14}. Let $R$ be a weakly left localizable ring such that the ring $R':= R/\gll_R$ satisfy the conditions of Theorem \ref{26Mar14}. Then
\begin{enumerate}
\item  the map $\phi : \maxDen_l(R)\ra \maxDen_l(R')$, $S\mapsto \pi'(S)$, is a bijection. In particular,
    $|\maxDen_l(R)|=|\maxDen_l(R')|<\infty$.
\item $\pi'^{-1}(\CC_{R'}) =\CC_l(R)$ and $\pi'(\CC_l(R)) = \CC_{R'}$.
\item If $\gll_R\subseteq \ass (\CC_l(R)):=\{ r\in R\, | \, cr=0$ for some $c\in \CC_l(R)\}$ then $\CC_l(R)\in \Den_l(R, \gll_R)$ and $ \CC_l(R)^{-1}R \simeq Q_{l, cl }(R)$.
\end{enumerate}
\end{proposition}

{\it Proof}. 1. By Theorem \ref{26Mar14}, $Q':= Q_{l, cl}(R') = \prod_{S'\in \maxDen_l(R')}S'^{-1}R'$ and ${}_RR'$ is an essential $R$-submodule of $Q'$. Since $R'\subseteq Q:= \prod_{S\in \maxDen_l(R)}S^{-1}R \subseteq Q'$ (Lemma \ref{a28Mar14}.(1)), $Q$ is an essential $R'$-submodule of $Q'$.  Then it follows from $$Q'= Q\times\prod_{S'\in \maxDen_l(R')\backslash \in \, \im \, \phi}S'^{-1}R'$$ that $\phi$ is a surjection, i.e. a bijection.

2. By Corollary  \ref{b26Mar14}.(2), $\CC_{R'}= \CC_l(R')$. Now, statement 2 follows from statement 1 and Lemma \ref{a28Mar14}.(3).

3. Since $ \CC_l(R)= \bigcap_{S\in \maxDen_l(R)} S$, $\ass (\CC_l(R))\subseteq \bigcap_{S\in \maxDen_l(R)} \ass (S) = \gll_R$. Hence, $\ass (\CC_l(R))= \gll_R$ (as $\ass (\CC_l(R))\supseteq \gll_R$, by the assumption). By statement 2, $\pi'(\CC_l(R))= \CC_{R'}$ and $\CC_{R'}\in \Den_l(R', 0)$. These two equalities together with $\ass (\CC_l(R))= \gll_R$ imply that $\CC_l(R)\in \Den_l(R, \gll_R)$. Then, clearly, $\CC_l(R)^{-1} R\simeq \pi'(\CC_l(R))^{-1} R'\simeq \CC_{R'}^{-1}R' = Q_{l,cl}(R)$.
$\Box $


{\bf Characterization of a class of weakly left localizable rings}. The next theorem is a characterization, via their left quotient rings,  of a class of weakly left localizable rings that satisfy some natural conditions.
\begin{theorem}\label{28Mar14}
Let $R$ be a ring, $\gll = \gll_R$, $\pi': R\ra R':= R/\gll$, $r\mapsto \br := r+\gll$.
The following statements are equivalent.
\begin{enumerate}
\item The ring $R$ is a weakly left localizable ring such that
\begin{enumerate}
\item the map $\phi : \maxDen_l(R)\ra  \maxDen_l(R')$, $S\mapsto \pi' (S)$, is a surjection.
\item $|\maxDen_l(R)|<\infty$,
\item for every $S\in \maxDen_l(R)$, $S^{-1}R$ is a weakly left localizable ring, and
    \item for all $S, T\in \maxDen_l(R)$ such that $S\neq T$, $\ass (S)$ is not a nil ideal modulo $\ass (T)$.
\end{enumerate}
\item $Q_{l, cl}(R')\simeq \prod_{i=1}^n R_i$ where $R_i$ are local rings with $\rad (R_i) = \CN_{R_i}$, $\gll$ is a nil ideal and $\pi'(\CL_l(R))= \CL_l(R')$.
\item $Q_l(R')\simeq \prod_{i=1}^n R_i$ where $R_i$ are local rings with $\rad (R_i) = \CN_{R_i}$,  $\gll$ is a nil ideal and $\pi'(\CL_l(R))= \CL_l(R')$.
\end{enumerate}
\end{theorem}

{\it Proof}. $(1\Leftrightarrow 2)$ By Proposition \ref{c13Dec12}, the ring $R$ is a weakly left localizable ring iff the ring $R'$ is a weakly left localizable ring, $\gll$ is a nil ideal and $\pi'(\CL_l(R))= \CL_l(R')$. Now, the equivalence $(1\Leftrightarrow 2)$ follows from Proposition \ref{b28Mar14}.(1), Theorem \ref{26Mar14} and Lemma \ref{a26Mar14}. Let give more details.

$(1\Rightarrow 2)$  {\em We have to show that the conditions (a)-(d) of Theorem \ref{26Mar14}.(1)  hold for the ring} $R'$. The condition (a) is obvious as $\gll_{R/ \gll} =0$ for any ring $R$ (Proposition \ref{a14Dec12}.(1)). The condition (b) of Theorem \ref{26Mar14} follows from the assumptions (a) and (b) of statement 1. By the assumption (a) of statement 1 and Lemma \ref{a28Mar14}.(1), the map $\phi$ is a bijection. The condition (c) of Theorem \ref{26Mar14} holds since, for all $S'\in \maxDen_l(R')$, $S'^{-1}R'\simeq S^{-1}R$  (where $S=\phi^{-1} (S')$) is a weakly left localizable ring, by the assumption (c) of statement 1. The condition (d) of Theorem \ref{26Mar14}  holds since $\phi$ is a bijection and, for all $S'\in \maxDen_l(R')$, $\ass (S') = \ass (S)/ \gll$ where $S= \phi^{-1}(S')$, $\ass (S)$ is not a nil ideal (by the assumption (d) of statement 1 and the fact that $\gll$ is a nil ideal.

$(2\Rightarrow 1)$   Since $Q_{l,cl}(R')=\prod_{i=1}^n R_i$ where $R_i$ are local rings with $\rad (R_i) = \CN_{R_i}$, by Theorem \ref{26Mar14}, the conditions (a)-(d) of Theorem \ref{26Mar14} hold and $R'$ is a weakly left localizable ring. Since $R'$ is a weakly left localizable ring, $\gll$ is a nil ideal and $\pi'(\CL_l(R)) =\CL_l(R')$, the ring $R$ is a weakly left localizable ring, by Proposition \ref{c13Dec12}. By Proposition \ref{b28Mar14}.(1), the map $\phi$ is a bijection, i.e. the condition (a) holds. Then the conditions (b) and (c) follow from the conditions (b) and (c) of Theorem \ref{26Mar14}.(1), respectively. Finally, the condition (d)  holds since the map $\phi$ is a bijection and for all $S\in \maxDen_l(R)$, $\ass (\phi (S)) = \ass (S)/ \gll$, $\ass (\phi (S))$ is not a nil ideal modulo $\ass (\phi (T))= \ass (T) / \gll$ for all $T\in \maxDen_l(R)$ such that $T\neq S$  (by the condition (d) of Theorem \ref{26Mar14}.(1) for the ring $R'$ and  since $\gll$ is a nil ideal).

$(2\Rightarrow 3)$ The implication is obvious as $Q_{l,cl}(R)=Q_l(R)$.

$(3\Rightarrow 2)$ The (left) $R$-module $R$ is an essential submodule of the (left) $R$-module $Q_l(R) = \prod_{i=1}^n R_i$. Therefore, for every $c= (c_1, \ldots , c_n) \in \CC_R$, where $c_i\in R_i$, the $R$-module homomorphism
$$ \cdot c: Q_l(R)\ra Q_l(R), \;\; q=(q_1, \ldots , q_n)\mapsto qc= (q_1c_1, \ldots , q_nc_n), $$ is an injection. Therefore, $c_i\in R_i^*$ for $i=1, \ldots , n$, i.e. $c\in Q_l(R)^*$. Therefore, $\CC_R\in \Den_l(R, 0)$ and $Q_l(R)= Q_{l, cl}(R)$. The proof of the theorem is complete.  $\Box $

$\noindent $

Corollaries \ref{c28Mar14} and \ref{d28Mar14} show that the rings that satisfy the conditions of Theorem \ref{28Mar14} have interesting properties and there is tight connection between rings $R$ and $R/\gll_R$.
\begin{corollary}\label{c28Mar14}
We keep  the notation of Theorem \ref{28Mar14}. Suppose that a ring $R$ satisfies one of the equivalent conditions 1--3 of Theorem \ref{28Mar14}. Then
\begin{enumerate}
\item $\maxDen_l(R) = \{ S_1, \ldots , S_n\}$ where $S_i = \{ r\in R\, | \, \frac{r}{1}\in R_i^*\}$, $S_i^{-1}R=R_i$ and the map $\phi$ is a bijection.
\item $\CC_l(R) = \pi'^{-1}(\CC_{R'})$ and $\pi'( \CC_l(R))= \CC_{R'}$.
\item $\pi'(\CL_l(R)) = \CL_l(R')$ and $\pi'^{-1} (\CL_l(R')) = \CL_l(R)$.
\item $\pi'(\Nil (R)) = \CN_{R'}$ and $\pi'^{-1} (\CN_{R'})= \Nil (R)$.
\item $\Nil (R) = \CN_R$.
\end{enumerate}
\end{corollary}

{\it Proof}. 1. In the proof of the implication $(1\Rightarrow 2)$ of Theorem \ref{28Mar14}, it was proven that the conditions (a)--(d) of Theorem \ref{26Mar14}.(1) hold for the ring $R' = R/ \gll$. By Proposition \ref{b28Mar14}.(1), the map $\phi$ is a bijection. Now,
 statement 1 follows from Theorem \ref{15Nov10} and Corollary \ref{b26Mar14}.(1).

 2. Suppose that one of the equivalent conditions 1--3 of Theorem \ref{28Mar14} holds. Then, by Theorem \ref{26Mar14},
\begin{eqnarray*}
\pi'^{-1}(\CC_{R'}) &=& \pi'^{-1} (\bigcap_{S'\in \maxDen_l(R')}S')=
 \bigcap_{S'\in \maxDen_l(R')}\pi'^{-1}(S) \\
 &=& \bigcap_{S\in \maxDen_l(R))}S=\CC_l(R),
\end{eqnarray*}
since $\phi$ is a bijection. The map $\pi'$ is a surjection. So, the first equality in statement 2 implies the second one.

3--5. The ring $R$ is a weakly left localizable ring. So, $R= \CL_l(R)\coprod \Nil (R)$. By Theorem \ref{28Mar14}.(2), Theorem \ref{26Mar14} and Corollary \ref{b26Mar14}.(3), $R'= \CL_l(R')\coprod \CN_{R'}$. By Theorem \ref{28Mar14}.(2), $\pi'(\CL_l(R))= \CL_l(R')$. Since $\CL_l(R)+\gll \subseteq \CL_l(R)$ (Proposition \ref{a14Dec12}.(2)), $\pi'^{-1}(\CL_l(R')) = \CL_l(R)$. Then $\pi'(\Nil (R))= \CN_{R'}$ and $\pi'^{-1} (\CN_{R'})= \Nil (R)$ since $\gll$ is a nil ideal of $R$. Since $\gll$ is a nil ideal of $R$, the ideal $\pi'^{-1} (\CN_{R'})$ is a nil ideal of $R$, i.e. $\Nil (R) = \CN_R$.  $\Box $


{\bf The set of completely left localizable elements is a left denominator set}. The next corollary gives sufficient conditions for the set $\CC_l(R)$ of completely left localizable elements of a ring $R$  that satisfies the conditions of Theorem \ref{28Mar14} to be a left denominator set.
\begin{corollary}\label{d28Mar14}
We keep of the notation of Theorem \ref{28Mar14}. Suppose that a ring $R$ satisfies one of the equivalent conditions 1--3 of Theorem \ref{28Mar14} and $\gll \subseteq \ass (\CC_l(R)):= \{ r\in R\, | \, cr=0$ for some $c\in \CC_l(R)\}$. Then
\begin{enumerate}
\item $\CC_l(R)\in \Den_l(R, \gll )$ and $\CC_l(R)^{-1}R\simeq Q_{l,cl}(R')$.
\item $\CC_l(R)^{-1}\CN_R\simeq \CC_{R'}^{-1} \CN_{R'}= \CN_{Q_{l,cl}(R')}= \rad (Q_{l,cl}(R'))$ and $\CN_R=\tau^{-1}(\CN_{Q_{l,cl}(R')})$ where $\tau : R\ra \CC_l(R)^{-1}R\simeq Q_{l,cl}(R')$, $r\mapsto \frac{r}{1}$ (see statement 1).
\item $\CC_l(R)^{-1}\CL_l(R)= \CL_l(\CC_l(R)^{-1}R)$ and $\tau^{-1}(\CL_l(Q_{l,cl}(R')))=\CL_l(R)$.
    \item $\CC_l(R)=S_{l,\gll }(R)$ and $Q_{l, \gll}(R)\simeq Q_{l,cl}(R')$.
        \item $S_l(R)^{-1}\gll$ is an ideal of $Q_l(R)$ and $Q_{l,cl}(Q_l(R)/S_l(R)^{-1}\gll )\simeq Q_{l,cl}(R')$.
\end{enumerate}
\end{corollary}

{\it Proof}. 1. Since $\CC_l(R) = \bigcap_{S\in \maxDen_l(R)} S$, we have the inclusion $$\ass (\CC_l(R))\subseteq \bigcap_{S\in \maxDen_l(R)} \ass (S)=\gll_R,$$ hence $\ass (\CC_l(R)) = \gll$ (since $\gll \subseteq \ass (\CC_l(R))$, by the assumption). By Corollary \ref{c28Mar14}.(2), $\pi'(\CC_l(R)) = \CC_{R'}\in \Den_l(R', 0)$ (Theorem \ref{28Mar14}.(2)). This inclusion together with the equality $\ass (\CC_l(R)) = \gll$ implies $\CC_l(R)\in \Den_l(R, \gll )$, Then, clearly,
$$\CC_l(R)^{-1}R\simeq (\pi'(\CC_l(R))^{-1}\pi'(R)\simeq \CC_{R'}^{-1}R'=
Q_{l,cl}(R').$$

2. By statement 1, $\CC_l(R)^{-1} \CN_R= \pi'(\CC_l(R))^{-1} \pi'(\CN_R) = \CC_{R'}^{-1} \CN_{R'}$, by Corollary \ref{c28Mar14}.(4,5). By Corollary \ref{b26Mar14}.(5), $\CC_{R'}^{-1} \CN_{R'}= \CN_{Q_{l,cl}(R')}= \rad (Q_{l,cl}(R'))$. Notice that $\tau = \s'\pi':R\stackrel{\pi'}{\ra}R'\stackrel{\s'}{\ra}\CC_l(R)^{-1}R\simeq Q_{l,cl}(R')$ (statement 1) where $\s'(r):=\frac{r}{1}$. Since $\gll$ is a nil  ideal of $R$ (Theorem \ref{28Mar14}), $\CN_{R'}=\pi'(\CN_R)$ and $\CN_R=\pi'^{-1}(\CN_{R'})$. Since $\CC_l(R)\in \Den_l(R, \gll )$ (statement 1) and $\gll \subseteq \CN_R$ (Theorem \ref{28Mar14}), we have $\CN_{R'}=\s'^{-1}(\CC_l(R)^{-1}\CN_R)= \s'^{-1}(\CN_{Q_{l,cl}(R')})$. Now,
$\CN_R=\pi'^{-1}(\CN_{R'})=(\s'\pi')^{-1}(\CN_{Q_{l,cl}(R')})=
\tau^{-1}(\CN_{Q_{l,cl}(R')})$.

3. By statement 1,
\begin{eqnarray*}
\CC_l(R)^{-1} \CL_l(R) &= & \pi'(\CC_l(R)) ^{-1} \pi'(\CL_l(R)) =\CC_{R'}^{-1} \CL_l(R')\\
 &=&  \CL_l(\CC_{R'}^{-1}R')\;\;\; ({\rm by\;\; Corollary}\; \ref{b26Mar14}.(6))\\
 &=& \CL_l(\CC_l(R)^{-1}R) \;\;\; ({\rm by\;\; statement }\; 1).\\
 \tau^{-1}(Q_{l,cl}(R'))&=&\pi'^{-1}\s'^{-1}(\CC_l(R)^{-1}\CL_l(R)) =\pi'^{-1}(\CL_l(R'))\\
 &=& \CL_l(R) \;\;\;\;\;\; ({\rm Corollary}\;\; \ref{c28Mar14}.(3)).
\end{eqnarray*}

4. By statement 1, $\CC_l(R)\subseteq S_{l, \gll}(R)$. Since $\ass (S_{l, \gll}(R))=\gll\subseteq \ass (S)$ for all $S\in \maxDen_l(R)$, the semigroup generated by $S_{l, \gll}(R)$ and $S$ is a left denominator set of $R$ necessarily equal to $S$. Hence, $S_{l, \gll}(R)\subseteq S$, and so $S_{l, \gll}(R)\subseteq \CC_l(R)$. Therefore, $S_{l, \gll}(R)=\CC_l(R)$. Hence, $Q_{l, \gll}(R)\simeq Q_{l,cl}(R')$, by statement 1.

5. By Proposition \ref{A8Dec12}.(1), $S_l(R)\subseteq \CC_l(R)$. The kernel of the ring homomorphism $Q_l(R)=S_l(R)^{-1}R\ra \CC_l(R)^{-1}R\simeq Q_{l,cl}(R')$, $s^{-1}r\mapsto s^{-1}r$, is equal to $S_l(R)^{-1}\gll$. So, $S_l(R)^{-1}\gll$ is an ideal of $Q_{l,cl}(R)$ such that $R'\subseteq Q_l(R)/S_l(R)^{-1}\gll\simeq \pi'(S_l(R))^{-1}R'\subseteq Q_{l,cl}(R')$. Hence, $Q_{l,cl}(Q_l(R)/S_l(R)^{-1}\gll )$ $\simeq Q_{l,cl}(R')$. $\Box$


\section{The cores of maximal left denominator sets of a weakly left localizable ring}\label{TCLOS}

The aim of this section is to find the cores  of maximal left denominator sets of a weakly left localizable ring. In \cite{Bav-LocArtRing}, several properties of  the {\em core} $S_c$ of a left Ore set $S$ of a ring $R$ were established.

$\noindent $

{\it Definition}, \cite{Crit-S-Simp-lQuot}. Let $R$ be a ring and $S$ be its left Ore set. The subset of $S$,
$$ S_c:=\{ s\in S\, | \, \ker (s\cdot ) = \ass (S)\},$$ is called the {\em core} of the left Ore set $S$ where $s\cdot : R\ra R$, $r\mapsto sr$.

$\noindent $


Recall that a ring $R$ is called a {\em local ring} if $R\backslash R^*$ is an ideal of $R$ (equivalently, $R/ \rad (R)$ is a division ring). The next lemma is a criterion for  the ring $S^{-1}R$ (where $S\in \maxDen_l(R)$) to be a local ring.

\begin{lemma}\label{a20Apr14}
Let $R$ be a ring and $S\in \maxDen_l(R)$. Then the ring $S^{-1}R$ is a local  ring iff $R\backslash S$ is an ideal of the ring $R$.
\end{lemma}


{\it Proof}. Let $\s : R\ra S^{-1}R$, $ r\mapsto\frac{r}{1}$, and
$(S^{-1}R)^*$ be the group of units of the ring $S^{-1}R$. By
Theorem \ref{15Nov10}.(3), $ S= \s^{-1} ((S^{-1}R)^*)$.

$(\Rightarrow )$  If $S^{-1}R$ is a local ring then $S^{-1}R = (S^{-1}R)^*\coprod \gr$ where $\gr := \rad (S^{-1}R)$. Then
\begin{equation}\label{RSsr}
 R= \s^{-1}(S^{-1}R) = \s^{-1}((S^{-1}R)^*)\coprod \s^{-1}(\gr ) = S\coprod \s^{-1}(\gr ),
\end{equation}
and so $R\backslash S = \s^{-1}(\gr )$ is an ideal of $R$.

$(\Leftarrow )$  Suppose that $\gb := R\backslash S$ is an ideal of the ring $R$. Then $S^{-1}R = S^{-1}S\cup S^{-1}\gb $ where $S^{-1} S:= \{ s^{-1} s'\, | \, s,s'\in S\} = (S^{-1}R)^* $ (Theorem \ref{15Nov10}.(5)) and $S^{-1}\gb$ is a left ideal of the ring $S^{-1}R$.

(i) $S^{-1}S\cap S^{-1}\gb =\emptyset$: Suppose that the intersection is a non-empty set, i.e. $s^{-1}s'= t^{-1}b\in  S^{-1}S\cap S^{-1}\gb$ for some elements $s,s',t\in S$ and $b\in \gb$, we seek a contradiction. Then
$$ R\backslash S\ni b \in \s^{-1}(ts^{-1}s')\in \s^{-1} ((S^{-1}R)^*)=S,$$
a contradiction.

(ii) {\em For all $s\in S$, $\gb s^{-1}\subseteq S^{-1}\gb$,  i.e. $S^{-1}\gb = S^{-1}R\backslash (S^{-1}R)^*$ is an ideal of $S^{-1}R$, i.e. $S^{-1}R$ is a local ring}: Suppose that $\gb s^{-1}\not\subseteq S^{-1}\gb$ for some $s\in S$, i.e. $bs^{-1}\not\in S^{-1}\gb$ for some element $b\in \gb$. By the statement (i), $bs^{-1}\in S^{-1}S= (S^{-1}R)^*$, and so
$R\backslash S \ni b\in \s^{-1}((S^{-1}R)^*) =S$, a contradiction.   $\Box $


The next theorem gives an explicit description of the cores of maximal left denominator sets of a weakly left localizable ring that satisfies the conditions of Theorem \ref{26Mar14}.

\begin{theorem}\label{C2Dec12}
Let $R$ be a  ring  that satisfies the conditions of Theorem \ref{26Mar14}. Let  $\maxDen_l(R) = \{ S_1, \ldots , S_n\}$.
\begin{enumerate}
\item If $n=1$ then $S_{1,c} = S_1= R\backslash \CN_R$.
\item If $n\geq 2$ then $S_{i, c} = S_i\cap \bigcap_{j\neq i} \ga_i\neq \emptyset $  where $\ga_j=\ass (S_j)$.
\end{enumerate}
\end{theorem}

{\it Proof}. We keep the notation of Theorem \ref{26Mar14} and its proof.

1. If $n=1$ then $R=S_1\coprod \Nil (R)$, $\Nil (R) = \CN_R$ (Corollary \ref{b26Mar14}.(3)), $S_1= R\backslash \CN_R\in \Den (R, 0)$
(Theorem \ref{26Mar14}), and so $S_{1, c}=  S_1$.

2. Suppose that $n\geq 2$. For each $i=1, \ldots , n$, $\CC_i':=
S_i\cap \bigcap_{j\neq i} \ga_j\neq \emptyset$, see the statement (ii) in the proof of Theorem \ref{26Mar14}.  By Theorem \ref{26Mar14}, the map
$$ \s : =\prod_{i=1}^n \s_i : R\ra \prod_{i=1}^n R_i , \;\; r\mapsto (r_1, \ldots , r_n), $$ is a ring monomorphism.
 Since $R_i=S_i^{-1}R$ is a local ring  with $\rad (R_i) = \CN_{R_i}$ (Theorem \ref{26Mar14}),  by Lemma \ref{a20Apr14} and (\ref{RSsr}),
 $$R=S_i\coprod \s^{-1}_i (\CN_{R_i}).$$
 Clearly, $\ga_i =\{ r= (r_1, \ldots , r_n)\in R \, | \, r_i=0\}$. Each element $s'$ of the set $\CC_i'$ has the form
 $(0, \ldots , 0 , s_i', 0, \ldots , 0)$ with $s_i'\in R_i^*$. Then clearly, $s'\ga_i=0$ and so $\CC_i'\subseteq S_{i, c}$.
  To show that the equality $S_i'= S_{i, c}$ holds it suffices to show that every element $s\in S_i\backslash \CC_i'$ does
  not not belong to $S_{i, c}$. Fix $s$ such that $s\in S_i\backslash \CC_i'$. Then there is an index, say  $j$,  such that $j\neq i$ and
  such that $s_j\neq 0$ in $s= (s_1, \ldots , s_n)$. Then $s\cdot \CC_j' \neq 0$ but $\CC_j'\subseteq \ga_i$. Therefore, the element $s$
  does not belong to $S_{i, c}$.  $\Box $


\begin{corollary}\label{a22Apr14}
Let $R$ be a ring that satisfies the condition of Theorem \ref{28Mar14} and $\pi' : R\ra R' := R/ \gll_R$, $r\mapsto r+\gll_R$. Then, for all $S\in \maxDen_l(R)$, $\pi'(S_c) \subseteq \pi'(S)_c$.
\end{corollary}

{\it Proof}. By Corollary \ref{c28Mar14}.(1), the map $\phi : \maxDen_l(R)\ra \maxDen_l(R')$, $S\mapsto \pi'(S)$, is a bijection. Moreover, $\ass (\pi'(S))= \ass (S) / \gll_R$. Hence, $\pi'(S_c)\subseteq \pi'(S)_c$.  $\Box $



\section{Criterion for a semilocal ring to be a  weakly left localizable ring}\label{TTTBB}

A ring $R$ is called a {\em semilocal ring} if $R/\rad (R)$ is a semisimple (Artinian) ring.

$\noindent $

The next theorem is a criterion for a semilocal ring $R$ to be a weakly left localizable ring with  $\rad (R) = \CN_R$.

\begin{theorem}\label{24Dec12}
Let $R$ be a semilocal ring. Then the ring $R$ is a weakly left localizable ring with $\rad (R)=\CN_R$  iff $R\simeq \prod_{i=1}^s R_i$ where $R_i$ are local rings with $\rad (R_i) =\CN_{R_i}$.
\end{theorem}

{\it Proof}. $(\Leftarrow )$ Suppose that $R=\prod_{i=1}^s R_i$ is a direct product of local rings with $\rad (R_i) =\CN_{R_i}$. Then $\rad (R) =\prod_{i=1}^s \rad (R_i)$ and $\CN_R=\prod_{i=1}^s \CN_{R_i}$. Then the  equality $\rad (R)= \CN_R$
  follows from the equalities $\rad (R_i) = \CN_{R_i}$ for $i=1, \ldots , s$.  The local  rings $R_i=R_i^*\coprod\rad (R_i)=R_i^*\coprod\CN_{R_i}=R_i^*\coprod \Nil (R)$ are weakly left localizable rings. Hence, so is their direct product $R$, by Theorem  \ref{9Feb13}.

$(\Rightarrow )$ The ring $R$ is a semilocal ring, i.e. $\bR := R/ \rad (R) = \prod_{i=1}^s \bR_i$ where $\bR_i$ are simple Artinian rings, i.e. $\bR_i=M_{n_i}(D_i)$ is a ring of $n_i\times n_i$ matrices with coefficients from a division ring $D_i$ for $i=1, \ldots , s$. Let $\bE_{pq}(i)$ be the matrix units of $M_{n_i}(D_i)$ where $p,q=1, \ldots , n_i$. Since $\rad (R) = \CN_R$ is a nil ideal, the decomposition of 1 in $\bR$ as a sum of orthogonal idempotents
$$1=\sum_{i=1}^s\sum_{j=1}^{n_i} \bE_{jj}(i),$$ can be lifted to a decomposition of $1$ as a sum of orthogonal idempotents in the ring $R$,
$$1=\sum_{i=1}^s\sum_{j=1}^{n_i} E_{jj}(i).$$
(i) $n_1=\cdots = n_s$: If not then, say $n_i\geq 2$, we seek a contradiction. The nonzero idempotent $e_1=E_{11}(1)$ is a left localizable element of the ring $R$, i.e. $e_1\in S$ for some $S\in \Den_l(R)$. Since $e_1E_{ii}(j)=0$ for all idempotents $E_{ii}(j)$ distinct from $e_1$, we must have $E_{ii}(j)\in \ga := \ass (S)$. Let $J$ be the ideal in $R$ generated by the idempotent $E_{ii}(j)$ such that $E_{ii}(j)\neq e_1$. Then $J\subseteq \ga$ and $J+\rad (R) = R$ (since $n_i\geq 2$ and the rings $R_j$ are simple rings), and so $j'+r=1$ for some elements $ j'\in J$ and $r\in \rad (R)$. This implies that the element $j'=1-r$ is a unit in $R$, and so $\ga = R$, a contradiction.

(ii) {\em The idempotents $E_{11}(i)$, $i=1, \ldots , s$ are central idempotents}: Notice that $1=E_{11}(1)+\cdots +E_{ss}(s)=e_1+e_2$ is the sum of orthogonal idempotents. The ring $R$ can be seen as a matrix ring
$$ R=\bigoplus_{i,j=1}^2R_{ij}=\begin{pmatrix}
 R_{11}& R_{12}\\ R_{21} & R_{22}
 \end{pmatrix}
 \;\; {\rm where}\;\; R_{ij}:=e_iRe_j.$$
Let $S$  and $\ga$ be as above. Since $e_1(R_{21}+R_{22})=0$ and $R_{12}e_1=0$, we must have the inclusion $R_{21}+R_{22}+R_{12}\subseteq \ga$ (as $e_1\in S$ and $S\in \Den_l(R))$. Since $ R_{11}$ is a local ring
  and $\rad (R_{11})$ is a nil ideal (since $\rad (R_{11})\subseteq \rad (R)$) we must have the inclusion
  $$ S\subseteq \begin{pmatrix}
 R_{11}^*& R_{12}\\ R_{21} & R_{22}
 \end{pmatrix}$$ where $R_{11}^*$ is the group  of units of the ring $R_{11}$.
  In more detail, suppose that the inclusion does not hold, that is there exists an element
  $s= \begin{pmatrix}
 \l & \mu \\ \nu & \d
 \end{pmatrix}\in S\cap \begin{pmatrix}
 \rad (R_{11})& R_{12}\\ R_{21} & R_{22}
 \end{pmatrix}$. Using the fact  (Corollary 4.3, \cite{Bav-LocArtRing}) that if $T\in \Den_l(R, \ass (T))$ then $T+\ass (T) \in \Den_l(R, \ass (T))$, we may replace $S$ by $S+\ga$ and then we may assume that $S+\ga \subseteq S$. Since $ \begin{pmatrix}
 0 & \mu \\ \nu & \d
 \end{pmatrix}\in R_{21}+R_{22}+R_{12}\subseteq \ga$, the element
 $s'= \begin{pmatrix}
 \l & 0 \\ 0 & 0
 \end{pmatrix}\in S\cap \begin{pmatrix}
 \rad (R_{11})& R_{12}\\ R_{21} & R_{22}
 \end{pmatrix}$. This is impossible as $\l \in \rad (R_{11})$ is a nilpotent element.

  We claim that $R_{12}=0$. Let $a\in R_{12}$, i.e. $a=\begin{pmatrix}
 0& a\\ 0 & 0
 \end{pmatrix}$ then $sa=0$ for some elements $s=\begin{pmatrix}
 u& x\\ y & z
 \end{pmatrix}\in S$ where $u\in R_{11}^*$. Now,
 $$0=sa=\begin{pmatrix}
 u& x\\ y & z
 \end{pmatrix}\begin{pmatrix}
 0& a\\ 0 & 0
 \end{pmatrix} = \begin{pmatrix}
 0& ua\\ 0 & ya
 \end{pmatrix} .$$
 Then $a=0$ since $u\in R_{11}^*$, i.e. $R_{12}=0$. Using the same argument but for the idempotent $e_2$, we obtain that $R_{21}=0$. This means that $e_1$ and $e_2$ are central orthogonal idempotents. By symmetry, the idempotents $E_{11}(1), \ldots , E_{11}(s)$ are central. This means that $R\simeq \prod_{i=1}^sR_{ii}$ is the direct product product of local rings. Since $\rad (R)=\CN_R$, we must have $\rad (R_i) = \CN_{R_i}$ for all $i$.  $\Box $


\begin{corollary}\label{a24Decd12}
Let $R$ be a left Artinian ring. Then $R$ is a weakly left localizable ring iff  $R=\prod_{i=1}^sR_i$ is a direct product of local left Artinian rings.
\end{corollary}

{\it Proof}. Every left Artinian ring is a semilocal ring with $\rad (R) = \CN_R$. Now, the corollary follows from Theorem \ref{24Dec12}.
 $\Box $


\begin{corollary}\label{b24Dec12}
Let $R$ be a weakly left localizable, semilocal  ring with $\rad (R) = \CN_R$, i.e. $R=\prod_{i=1}^sR_i$ is a direct product of local rings, by Theorem \ref{24Dec12}. Let $1=e_1+\cdots + e_s$ be the corresponding sum of central orthogonal idempotents. Then
\begin{enumerate}
\item $\maxDen_l(R)=\{ S_1, \ldots , S_s\}$ where $S_i= R_1\times \cdots \times R_i^*\times\cdots \times R_s$ for $i=1, \ldots , s$; $\ass (S_i) = R_1\times \cdots \times 0\times\cdots \times R_s$ and $S_i^{-1}R\simeq R_i$.
\item The core of $S_i$ is $0\times \cdots \times 0\times  R_i^*\times 0 \times \cdots \times 0$.
\item For each $i=1, \ldots , s$, $S_{e_i}:=\{ 1, e_i\}\in \Den_l(R, \ass (S_i))$ and $S_{e_i}^{-1} R\simeq S_i^{-1}R\simeq R_i$.
    \item $\Nil (R) = \CN_R$.
\end{enumerate}
\end{corollary}

{\it Proof}. 1 and 2. The ring $R_i$ are local rings with $\rad (R_i)=\CN_{R_i}$. Therefore, the rings $R_i$ are left localization maximal rings, i.e. $\maxDen_l(R_i) = \{ R_i^*\}$. Now, statements 1 and 2 follow from Theorem \ref{c26Dec12}.





3 and 4. Trivial . $\Box $


$${\bf Acknowledgements}$$

 The work is partly supported by  the Royal Society  and EPSRC.


\small{

Department of Pure Mathematics

University of Sheffield

Hicks Building

Sheffield S3 7RH

UK

email: v.bavula@sheffield.ac.uk}

\end{document}